\documentclass[a4paper,9pt]{article}

\usepackage{amsmath,amscd,amssymb}

\def\h{\mathfrak{h}}
\def\n{\mathfrak{n}}
\def\g{\mathfrak{g}}
\def\b{\mathfrak{b}}
\def\C{\mathbb{C}}
\def\Z{\mathbb{Z}}

\newtheorem{theo}{\bf{Theorem}}

\newtheorem{lem}{Lemma}

\begin{document}
\title{On highest weight modules\\ over elliptic quantum groups}
\author{Pavel Etingof and Olivier Schiffmann\\
\texttt{etingof@math.harvard.edu}, \texttt{schiffma@clipper.ens.fr}}
\maketitle
The purpose of this note is to define and construct highest weight modules
for Felder's elliptic quantum groups. This is done by using exchange matrices 
for intertwining operators between modules over quantum affine algebras. 

A similar problem for the elliptic quantum group corresponding to 
Belavin's R-matrix was posed in \cite{FIJKMY}. This problem, as well as 
its analogue for Felder's R-matrix was solved in a recent paper 
\cite{JKOS}. Thus our goal is to suggest another solution of this problem. 
Our approach is similar to that of \cite{JKOS}
but somewhat different. Namely, we construct the quasi-Hopf twist
(which is necessary to pass from the quantum affine algebra 
to elliptic algebra) not as an infinite product of R-matrices, 
but axiomatically as a 
fusion matrix for intertwining operators (this forces us to consider
intetwining operators taking values in arbitrary, not necessarily 
finite dimensional modules).
For finite dimensional Lie algebras and quantum groups, 
a similar construction was introduced in \cite{EV1}. 
 The equivalence of the two approaches (ours and that of \cite{JKOS})
follows from the fact that the fusion matrix satisfies a version 
of the quantum KZ equations, which implies that it is a one-sided infinite 
product of (modified) R-matrices.    
 
\section{Quantum affine algebras and exchange dynamical quantum groups}
\subsection{The quantum affine algebra $U_q(\hat{\mathfrak{sl}}_n)$} 
\paragraph{}Consider $\g=\mathfrak{sl}_n(\C)$ equipped with the invariant form $(a,b)=tr(ab)$ and let $\h \subset \g$ be the Cartan subalgebra of diagonal matrices. Let $\g=\n_+ \oplus \h \oplus \n_-$ denote the usual polarization and let $(e_i),(f_i)$, $i=1,\ldots n-1$ be corresponding Chevalley generators. Set $\check{h}=n$ (the dual Coxeter number) and let $\rho \in \h^*$ be the half-sum of positive roots.
Let $\hat{\g}$ be the affine algebra associated to $\g$; we will denote its Chevalley generators by $(e_i),(f_i)$ for $i=0,\ldots n-1$, and by $c$ the central element (see \cite{Kac}).\\
\hbox to1em{\hfill}Fix some $q \in \C^*$ (not a root of unity) and 
a value of $\log q$. By $q^a$ we always mean $e^{a\log q}$. Let $U:=U_q(\hat{\g})$ denote the Drinfeld-Jimbo quantum affine (Hopf) algebra, generated by $(E_i),(F_i)$, for $i=0,\ldots n-1$, $K_i=q^{h_i}$ for $i=1,\ldots n-1$ and the central element $q^c$, with the quantum Serre presentation (as in \cite{CP}, chapter 6). We let $U^+:=U_q(\hat{\b}_+)$ (resp. $U^-:=U_q(\hat{\b}_-)$) be the Hopf subalgebra of $U$ generated by $E_i,K_j,q^c$ (resp. by $F_i,K_j,q^c$).\\
\hbox to1em{\hfill}It is useful to consider the extended quantum affine algebra $\tilde{U}:=U_q(\tilde{\frak g})$ generated by $U$ and an element $d$ 
 satisfying the relations $[d,E_j]=[d,F_j]=[d,K_j]=0$ for $j>0$, $[d,E_0]=E_0,\;[d,F_0]=-F_0$, $[d,K_0]=0$ and $\Delta(d)=d\otimes 1 + 1 \otimes d$, $\epsilon(d)=0$, $S(d)=-d$. We let $\tilde{U}^+$ (resp. $\tilde{U}^-$) denote the Hopf subalgebra of $\tilde U$ generated by $d$ and $U^+$ (resp $d$ and $U^-$). The Hopf algebra $\tilde{U}$ is (topologically) quasitriangular, and we will denote by $\mathcal{R}$ its universal R-matrix. It is known that $\mathcal{R} \in q^{d \otimes c + c\otimes d}U^+\hat{\otimes}U^-$
(here as usual $\hat\otimes$ denotes the completed tensor product).\\
\hbox to1em{\hfill}For any $\lambda \in \h^*$ and $k \in \C$, define the Verma module of highest weight $\lambda$ and level $k$ over $\tilde U$ by $M_{\lambda,k}=\mathrm{Ind}_{\tilde{U}^+}^{\tilde{U}}(\chi_{\lambda,k})$ where $\chi_{\lambda,k}$ is the one-dimensional $\tilde{U}^+$ module generated by a vector $x_{\lambda,k}$ such that $E_i.x_{\lambda,k}=0$ for $i\geq 0$, $K_i.x_{\lambda,k}=q^{\lambda(h_i)}x_{\lambda,k}$, $q^c.x_{\lambda,k}=q^kx_{\lambda,k}$ and $d.x_{\lambda,k}=-\Delta_k(\lambda)x_{\lambda,k}$ where we set $\Delta_k(\lambda)=\frac{(\lambda,\lambda + 2\rho)}{2(k+\check{h})}$ (see \cite{EFK}). We will denote by $M^*_{\lambda,k}$ the restricted dual Verma module and by $x_{\mu,k}^*$ its lowest weight vector such that $x^*_{\mu,k}(x_{\mu,k})=1$.\\
\hbox to1em{\hfill}Let $\mathrm{Rep}_f(U)$ denote the category of finite-dimensional $U$-modules $V$ with weight decomposition $V=\bigoplus_\mu V[\mu]$ where $q^{h_i}_{|V[\mu]}=q^{\mu(h_i)}$ and 
$q^c$ acts by 1. For $V \in \mathrm{Rep}_f(U)$,  and $z \in \C^*$, we denote by $V(z)$ the twist of $V$ by the automorphism $D_z:U \to U$ defined by $D_z(E_0)=zE_0$, $D_z(F_0)=z^{-1}F_0$ and $D=\mathrm{Id}$ on all other generators (i.e, we set $\pi_{V(z)}(a)=\pi_V(D_z(a))$ for $a \in U$). Similarly, for any $\Delta \in \C$, we let $z^{-\Delta}V[z,z^{-1}]:=V \otimes z^{-\Delta}\C[z,z^{-1}]$ be the $\tilde{U}$-module defined as follows:.
$U$ acts in this module by $\pi_{z^{-\Delta}V[z,z^{-1}]}(a)=
\pi_V(D_z(a))$, and $d$ acts as $z\frac{d}{dz}$
(of course, here $z$ is no longer a number but a formal variable).
\paragraph{}The action of the universal R-matrix $\mathcal{R} \in \tilde{U} \hat{\otimes} \tilde{U}$ on objects of $\mathrm{Rep}_f(U)$ is given by a power series
$$\mathcal{R}_{V, W}(\frac{z_1}{z_2})=(\pi_{V(z_1)} \otimes \pi_{W(z_2)})(\mathcal{R})$$
where $\pi_X$ denotes the action of $U$ on the module $X$. (Note that 
although $d$ is not defined on finite-dimensional representations, 
$\mathcal R$ is still defined, as $c=0$.)
This power series converges in the region $z_1\ll z_2$, i.e in the neighborhood of $\frac{z_1}{z_2}=0$, is regular at $0$, and admits a meromorphic continuation to $\frac{z_1}{z_2} \in \C$ (see \cite{KS}, Corollary 4.4.12).
\subsection{Intertwining operators and quantum correlation functions}
\paragraph{}Fix complex numbers $k,l\ne -\check h$ and let $Y$ be a
$\tilde U$-module in which the elements $K_i$, $q^c$, $d$ act diagonally.
We denote by $Y[\nu,r,s]$ the eigenspace of weight $\nu$ where $q^c$ acts
as $q^r$ and $d$ as $s$. 
The following theorem is a special case of 
Frobenius reciprocity (see e.g.\cite{FR}). 
\begin{theo}Let $\lambda,\mu \in \h^*$ be such that $M^*_{\mu,l}$ is irreducible. Then
$$\mathrm{Hom}_{\tilde U}(M_{\lambda,k},M_{\mu,l}\hat{\otimes} Y)\simeq Y[\lambda-\mu,k-l,\Delta_l(\mu)-\Delta_k(\lambda)]$$ where the isomorphism is as follows: to $\Phi: M_{\lambda,k} \to M_{\mu,k} \hat{\otimes} Y$ we associate the element \linebreak $<\Phi>=<x_{\mu,k}^*,\Phi(x_{\lambda,k})>\in Y$.
\end{theo}

{\bf Remark.} Here $M_{\mu,k}\hat\otimes Y$ denotes the tensor product
completed with respect to the grading of $M_{\mu,k}$ by affine weights. 

\paragraph{}The theorem allows us to assign to any vector 
$v \in Y[\lambda-\mu,k-l,\Delta_l(\mu)-\Delta_k(\lambda)]$
the corresponding intertwiner, which we denote by 
$\Phi_{\lambda,k}^v$. 

Now we will apply the theorem to 
the module $Y=z^{-\Delta}V[z,z^{-1}]$, $V\in  \mathrm{Rep}_f(U)$.  
The theorem imples that any vector $v\in V$ 
of weight $\nu$
there exists a unique 
intertwining operator 
$\Phi^{z^{-\Delta} v}_{\lambda,k}:
M_{\lambda,k} \to 
M_{\mu,k}\hat\otimes z^{-\Delta}V[z,z^{-1}]$,
where $\Delta=\Delta_k(\lambda)-\Delta_k(\mu)$, 
such that $<\Phi^{z^{-\Delta} v}_{\lambda,k}>=
z^{-\Delta}v$. For brevity we will denote this operator 
by $\tilde\Phi^v_{\lambda,k}(z)$. It can be represented as a series 
$\tilde\Phi^v_{\lambda,k}(z)=
\sum_{m \in \Z} \Phi_{\lambda,k}^v[m]z^{-m-\Delta}$, where
$\Phi_{\lambda,k}^v[n]$ is homogeneous of degree $n$. 

\paragraph{}Now let $V,W \in \mathrm{Rep}_f(U)$, 
$v \in V[\mu],\,w \in W[\nu]$ and consider the formal composition
$$\tilde{\Phi}_{\lambda-\nu,k}^v(z_1)\tilde{\Phi}_{\lambda,k}^w(z_2):M_{\lambda,k} \to M_{\lambda-\mu-\nu,k}\hat{\otimes}z_1^{-\Delta_1}V[z_1,z_1^{-1}] \hat{\otimes} z_2^{-\Delta_2}W[z_2,z_2^{-1}]$$
where $\Delta_1=\Delta_k(\lambda-\nu)-\Delta_k(\lambda-\nu-\mu)$; $\Delta_2=\Delta_k(\lambda)-\Delta_k(\lambda-\nu)$
(here $V[z_1,z_1^{-1}]\hat\otimes W[z_2,z_2^{-1}]$ is the tensor 
product completed with respect to the grading in the first component). 
This composition can be written as a series $\tilde{\Phi}_{\lambda-\nu,k}^v(z_1)\tilde{\Phi}_{\lambda,k}^w(z_2)=z_1^{-\Delta_1}z_2^{-\Delta_2}\sum_{m \in \Z_+} \Phi_{\lambda,k,m}^{v,w}(\frac{z_2}{z_1})^{m}$. Define the quantum correlation function as
$$\Psi_{\lambda,k}^{v,w}(z_1,z_2)=<x_{\lambda-\mu-\nu}^*,\tilde{\Phi}_{\lambda-\nu,k}^v(z_1)\tilde{\Phi}_{\lambda,k}^w(z_2).x_{\lambda}>.$$
Set $p=q^{-2(k + \check{h})}.$
\begin{theo}[\cite{FR}, see also \cite{EFK}] 
The formal power series $\Psi_{\lambda,k}^{v,w}(z_1,z_2)$ satisfies the qKZ equations
\begin{align}
\Psi_{\lambda,k}^{v,w}(z_1,pz_2)&=\mathcal{R}_{W,V}^{21}(\frac{pz_2}{z_1})q^{2\lambda-\nu-\mu+2\rho}_{|W}\Psi_{\lambda,k}^{v,w}(z_1,z_2)\label{E:qKZ1}\\
\Psi_{\lambda,k}^{v,w}(pz_1,z_2)&=q^{2\lambda-\nu-\mu+2\rho}_{|V}\mathcal{R}_{W,V}^{21}(\frac{z_2}{z_1})^{-1}\Psi_{\lambda,k}^{v,w}(z_1,z_2)\label{E:qKZ2}
\end{align}
Moreover,
\begin{equation}\label{E:01}
\Psi_{\lambda,k}^{v,w}(z_1,z_2)=z_1^{-\Delta_1}z_2^{-\Delta_2}F_{\lambda,k}^{v,w}(\frac{z_2}{z_1})
\end{equation}
where $F_{\lambda,k}^{v,w}(z)$ is a meromorphic function on $\C$
which is regular at the origin.
\end{theo}
\subsection{qKZ twist and dynamical quantum groups}

\paragraph{}Let us fix $k \in \C$ and let $V,W \in \mathrm{Rep}_f(U)$. Consider the map 
\begin{align*}
J_{V,W}(u_1,u_2,\lambda,k):\; V \otimes W & \to V \otimes W\\
v \otimes w &\mapsto 
e^{-2\pi i(\Delta_1u_1+\Delta_2u_2)}F_{\lambda,k}^{v,w}(e^{2\pi i(u_2-u_1)})
\end{align*}
This map depends meromorphically on $u_1,u_2$. 
We will call this map the fusion matrix, or the qKZ twist, since by Theorem 2 it satisfies
the quantum KZ equations.  
Now define the exchange matrix
$$R_{V,W}(u_1,u_2,\lambda,k)=J_{V,W}(u_1,u_2,\lambda,k)^{-1}\mathcal{R}^{21}_{W,V}(e^{2\pi i(u_2-u_1)})J^{21}_{W,V}(u_2,u_1,\lambda,k)$$
Notice that $R_{V,W}(u_1,u_2,\lambda,k)$ is of weight zero. It is easy to see using (\ref{E:01}) that $R_{V,W}(u_1,u_2,\lambda,k)$ depends only on 
$u=u_1-u_2$, so we will denote it by $R_{V,W}(u,\lambda,k)$. 
It is a meromorphic function of $u$.

\begin{theo} The operator $ R_{V,W}(u,\lambda,k)$ satisfies the quantum dynamical Yang-Baxter equation with a spectral parameter (see \cite{F}): for any $V,W,X \in \mathrm{Rep}_f(U)$, and $u_1,u_2,u_3 \in \C$ we have
\begin{equation}
\begin{split}
 R_{V,W}^{12}(u_1-u_2,\lambda-h^{(3)},k)& R^{13}_{V,X}(u_1-u_3,\lambda,k) R^{23}_{W,X}(u_2-u_3,\lambda-h^{(1)},k)\\
&= R_{W,X}^{23}(u_2-u_3,\lambda,k) R^{13}_{V,X}(u_1-u_3,\lambda-h^{(2)},k) R^{12}_{V,W}(u_1-u_2,\lambda,k)
\end{split}
\end{equation}
\end{theo}

This theorem is an analogue of the corresponding result for finite dimensional 
Lie algebras, which is proved in \cite{EV1}. The proof
is straightforward, as in the finite dimensional case. 

\paragraph{}We will set $  R_k(u,\lambda)=
 R_{\C^n,\C^n}(u,\lambda,k)$, where $\C^n$ denotes the standard vector representation. 
Theorem 3 implies that $ R_k(u,\lambda)$ is a
quantum dynamical R-matrix with a spectral parameter.
Moreover, it satisfies the modified unitarity condition:  
$ R_k(u,\lambda) R_k^{21}(-u,\lambda)=\chi(u)$,
where $\chi$ is an elliptic function with periods $1$ and $2i\pi^{-1}\log q$,
which does not depend on $\lambda$ and $k$. This follows from the fact that
$\mathcal{R}_{\C^n,\C^n}(z)$ 
satisfies a similar condition: 
$\mathcal R_{\C^n,\C^n}(z)\mathcal R_{\C^n,\C^n}^{21}(z^{-1})=
\chi(\frac{1}{2\pi i}\log z)$. 
(The function $\chi$ is computed in \cite{EFK},
\cite{FR}, but its
explicit form is not important for us).

\begin{lem} The function $ R_k(u,\lambda)$ satisfies the following 
periodicity properties:
\begin{enumerate}
\item $ R_k(u-\frac{i}{\pi}(k+\check h)\log q,\lambda)=\chi(u)^{-1}
 R_k(u,\lambda)$,
\item $ R_k(u+1,\lambda)=e^{2\pi i(\Delta_k(\lambda-h^{(2)})-\Delta_k(\lambda-h^{(1)}-h^{(2)}))} R_k(u,\lambda)e^{-2\pi i(\Delta_k(\lambda)-\Delta_k(\lambda-h^{(1)}))}$.
\end{enumerate}
\end{lem}
\textbf{Proof:} Relation 1. follows from the fact that $J_{\C^n,\C^n}$ 
satisfies the qKZ equations; relation 2. follows from (\ref{E:01}).$\square$

Thus, we see that $ R_k(u,\lambda)=\hat R_k(u,\lambda)\psi_k(u)$, where 
$\hat R_k$ is an ``elliptic'' unitary solution of the quantum Yang-Baxter
equation (i.e. a solution 
satisfying part 1 of Lemma 1 without the factor $\chi$, 
and part 2 of Lemma 1), and $\psi_k$ is a scalar meromorphic function. 
The precise form of $R_k$ will be determined in a forthcoming paper. 
In fact, it can be shown that $\hat R_k(u,\lambda)$ can be transformed into 
Felder's elliptic dynamical R-matrix defined in \cite{F} by gauge 
transformations of \cite{EV0}.

Denote by $\cal C$ the category of meromorphic representations 
of the R-matrix $R_k$, i.e. of pairs $(V,L_V)$, where 
$V$ is a diagonalizable finite dimensional $\h$-module, and 
$L_V:{\frak h}\times \C\to \text{End}(\C^n\otimes V)$ an invertible 
meromorphic function of weight zero, such that 
\begin{equation}
\begin{split}
 R_k^{12}(u_1-u_2,\lambda-h^{(3)})& L_V^{13}(u_1-u_3,\lambda) L_V^{23}(u_2-u_3,\lambda-h^{(1)})\\
&= L_V^{23}(u_2-u_3,\lambda) L_V^{13}(u_1-u_3,\lambda-h^{(2)}) R_k^{12}(u_1-u_2,\lambda)
\end{split}
\end{equation}
  
{\bf Remark.} This category does not change when $R_k$ is replaced 
with the "purely" elliptic R-matrix $\hat R_k$.
 
\begin{theo} The assignment 
$$
V  \mapsto  (V,L_V(u,\lambda)=R_{\C^n,V}(u,\lambda))
$$
and the identity at the level of morphisms defines a functor
$F:\; \mathrm{Rep}_f(U) \to \mathcal{C}$.
\end{theo}
\paragraph{Proof:} straightforward (see \cite{EV1}).

\subsection{Highest weight modules over $U_q(\hat{\mathfrak{sl}}_n)$}
\paragraph{}Consider the category ${\mathcal O}^l$ of 
$\tilde U$-modules of level $l \in \C$ from category $\mathcal{O}$. 
As usual, we denote by $W=\bigoplus_{\mu}W[\mu,s]$ the weight 
decomposition of such a module 
($\mu$ is the weight and $s$ the eigenvalue of $d$). 
For a complex number $s$, we denote by $W^s$ the module obtained from $W$ by 
shifting the action of $d$ by $s$. 

\paragraph{}The universal R-matrix $\mathcal{R} \in \tilde{U} \otimes \tilde{U}$ is defined on products of a module of the form 
$z^{-\Delta}V[z,z^{-1}]$, $V\in \mathrm{Rep}_f(U)$ with a module $W\in\mathcal O^l$:
$$\mathcal R_{V,W}(z)=(\pi_{z^{-\Delta}V[z,z^{-1}]} \otimes \pi_{W})(\mathcal{R}).$$ 
This is a Laurent series, which is independent of $\Delta$ up to scaling. 
Similarly, one defines 
$\mathcal R_{W,V}(z^{-1})=(\pi_{W}\otimes \pi_{z^{-\Delta}V[z,z^{-1}]})(\mathcal{R})$.  

Let $W$ be as above. By Theorem 1, for a generic $k$ and 
any homogeneous vector $w\in W[\nu,s]$, 
we can define 
an intertwining operator $\Phi^w_{\lambda,k}:
M_{\lambda,k}\to M_{\lambda-\nu,k-l}\otimes W^{\Delta_{k-l}(\lambda-\nu)-\Delta_l(\lambda)-s}$. 

\hbox to1em{\hfill}Now define the quantum correlation functions and the exchange matrices. 
 Let $V \in \mathrm{Rep}_f(U), v \in V[\mu]$, $W \in {\mathcal O}^l, w \in W[\nu]$ and consider the formal composition
$$\tilde{\Phi}_{\lambda-\nu,k-l}^v(z){\Phi}_{\lambda,k}^w:M_{\lambda,k} \to M_{\lambda-\mu-\nu,k-l}\hat{\otimes}z^{-\Delta}V[z,z^{-1}] \hat{\otimes}W^{\Delta_{k-l}(\lambda-\nu)-\Delta_l(\lambda)-s}$$
where $\Delta=\Delta_{k-l}(\lambda-\nu)-\Delta_{k-l}(\lambda-\nu-\mu)$. The quantum correlation function is
$$\Psi_{\lambda,k}^{v,w}(z)=<x_{\lambda-\mu-\nu,k-l}^*,\tilde{\Phi}_{\lambda-\nu,k-l}^v(z){\Phi}_{\lambda,k}^w.x_{\lambda,k}>.$$
\paragraph{}Define the qKZ twist
 \begin{align*}
J_{V,W}(z,\lambda,k):\; V \otimes W & \to V \otimes W\\
v \otimes w &\mapsto \Psi_{\lambda,k}^{v,w}(z)
\end{align*}
\hbox to1em{\hfill}Notice that this function is given by a series $J_{V,W}(z,\lambda,k)=z^{-\Delta}\sum_{n\in \Z} J_{V,W}[m]z^m$ which is infinite in both the positive and the negative directions (here $\Delta$ should be considered as a diagonal operator whose eigenvalues are given above). However it becomes finite in the negative direction when applied to an element $v \otimes w$.
\paragraph{} Using the same construction 
with the reversed order of tensor product, 
one defines the twist $J^*_{W,V}(z,\lambda,k):W\otimes V\to W\otimes V$. 
This twist is also a series infinite in both directions, but it becomes 
finite in the positive direction after being applied to a vector.  

\paragraph{}The exchange matrix is
$$\tilde R_{V,W}(z,\lambda,k)=J_{V,W}(z,\lambda,k)^{-1}
\mathcal{R}^{21}_{W,V}(z^{-1})J^{*21}_{W,V}(z,\lambda,k).$$
 It is a series of weight zero which is infinite in both directions, but 
becomes finite in the positive direction after being applied to a vector.  
Define $R_{V,W}(u,\lambda,k):=\tilde R(e^{2\pi iu},\lambda,k)$. 
\begin{theo}\label{Theo} The operator $R_{V,W}(u,\lambda,k)$ satisfies the quantum dynamical Yang-Baxter equation with central charge: for any $V,W \in \mathrm{Rep}_f(U)$, $X \in {\mathcal O}^l$, and $u,w \in \C$ we have
\begin{equation}
\begin{split}
R_{V,W}^{12}(u-u',\lambda&-h^{(3)},k-l)R^{13}_{V,X}(u,\lambda,k)R^{23}_{W,X}(u',\lambda-h^{(1)},k)\\
&=R_{W,X}^{23}(u',\lambda,k)R^{13}_{V,X}(u,\lambda-h^{(2)},k)R^{12}_{V,W}(u-u',\lambda,k)
\end{split}
\end{equation}
as operators $V \otimes W \otimes X \to V \otimes W \otimes X$.
\end{theo}

{\bf Remark.} This equality should be understood as one usually understands 
commutation relations between vertex operators: 
the two sides are series converging in two different regions
to the same meromorphic functions. 

The proof is analogous to the proof of Theorem 3.

\subsection{Modules over an elliptic quantum group 
with nontrivial central charge}
\paragraph{}If we take $V=W=\C^n$ (the vector representation) in Theorem~\ref{Theo}, we obtain the relation
\begin{equation}
\begin{split}
R_{k-l}^{12}(u-u',\lambda-h^{(3)})&R^{13}_{\C^n,X}(u,\lambda,k)R^{23}_{\C^n,X}(u',\lambda-h^{(1)},k)\\
&=R_{\C^n,X}^{23}(u',\lambda,k)R^{13}_{\C^n,X}(u,\lambda-h^{(2)},k)
R^{12}_{k}(u-u',\lambda)
\end{split}
\end{equation}
for any $X \in {\mathcal O}^l$. 
This motivates the following definition.

\paragraph{}Let ${\frak h}$ be a finite dimensional abelian Lie algebra, 
$\hat {\frak h}={\frak h}\oplus \C c$,
Elements of ${\frak{\hat h}}^*$ will be denoted by 
$\hat\lambda=(\lambda,k)$.  
Let $R(u,\hat\lambda)$ be a holomorphic 
family (parametrized by $k,\text{Re}\,k>s$) of 
quantum dynamical R-matrices on $\frak h^*$ 
with step 1 and spectral parameter $u$, 
with values in $\text{End}(\C^n\otimes \C^n)$, (see \cite{EV0}). 

\paragraph{Definition:} a {\it bounded (from above) representation of $R$ with central charge $l$} is a pair $(V,L(z,\hat\lambda))$ where $V=\bigoplus_{i<N} V[i]$ is a $\Z$-graded $\C$-vector space bounded from above, equipped with a diagonalizable $\hat\h$-action such that $c$ acts by $l$
(with finite dimensional homogeneous subspaces), 
and $L=(L_{ij})$, 
where 
$$L_{ij}(z,\hat\lambda)=\sum_{m \in \Z} L_{ij}(\lambda,k)[m]z^{-m-\Delta_{ij}}\in z^{-\Delta_{ij}}\mathrm{End}(\C^n \otimes V)[[z,z^{-1}]],\ (\Delta_{ij}\in\C)$$
is a Laurent series satisfying the following properties:
\begin{enumerate}
\item $L(\lambda,k)[n]$ is a weight zero homogeneous operator of degree $n$ depending meromorphically on $\lambda \in \h^*$ and $k$ in the region 
$\text{Re}\,k>s$.
\item for any $w \in \C^n$ and $v \in V$, $L(\lambda,k)[i](w \otimes v)=0$ for $i \ll 0$;
\item we have
\begin{equation*}
\begin{split}
R^{12}(u-u',\hat\lambda-\hat h^{(3)})&L^{13}(e^{2\pi iu},\hat\lambda)
L^{23}(e^{2\pi iu'},\hat\lambda-\hat h^{(1)})\\
&=L^{23}(e^{2\pi iu'},\hat\lambda)L^{13}(e^{2\pi iu},\hat\lambda-\hat h^{(2)})R^{12}(u-u',\hat\lambda)
\end{split}
\end{equation*}
i.e. the matrix elements of both sides coincide 
as meromorphic functions of $u,u' \in \C$, $\lambda \in \h^*$ and $k$.
\end{enumerate}

Some of the most interesting examples are:

1) ${\frak h}=0$, $R(u,k)$ is the Belavin R-matrix with elliptic modulus 
$\tau=ak+b$. 

2) ${\frak h}$ is like in Section 1, $R(u,\hat\lambda)$ is the Felder 
elliptic dynamical R-matrix with elliptic parameter $\tau=ak+b$. 

3) ${\frak h}$ is as in Section 1, and $R(u,\hat\lambda)=R_k(u,\lambda)$.   

As was mentioned above, one can show that examples 2 and 3 are essentially 
the same example, since the corresponding R-matrices are related by a gauge 
transformation.  

\paragraph{}A morphism between two modules $(V,L_V(u,\hat\lambda))$ and $(W,L_W(u,\hat\lambda))$ is a $\Z$-graded, $\h$-linear map $\varphi(\hat\lambda):V \to W$ such that
$$L_W(u,\hat\lambda)\varphi_2(\hat\lambda-\hat h^{(1)})=\varphi_2(\hat\lambda)L_V(u,\hat\lambda).$$
\paragraph{}Let us denote by $\mathcal{C}^l$ the category of bounded representations of the R-matrix $R(u,\hat\lambda):=R_k(u,\lambda)$ of central charge $l$. There is a notion of a tensor product between 
objects of $\mathcal C^{l}$ and $\mathcal C^{l'}$ which gives an object 
in $\mathcal C^{l+l'}$. Namely, if $(V,L_V(u,\hat\lambda))\in \mathcal{O}b(\mathcal{C}^l)$ and $(W,L_W(u,\hat\lambda))\in \mathcal{O}b(\mathcal{C}^{l'})$ then
$$(V \otimes W, L^{12}_V(u,\lambda-h^{(3)},k-l')L^{13}_W(u,\lambda,k))$$
is an object of $\mathcal{C}^{l+l'}$. This defines a bifunctor $\otimes: \mathcal{C}^l \times \mathcal{C}^{l'} \to \mathcal{C}^{l+l'}$ (with the tensor 
product of morphisms defined as in \cite{EV0}, i.e. 
$f\bar\otimes g(\hat\lambda)=f(\hat\lambda-\hat h^{(2)})\otimes 
g(\hat\lambda)$). The category $\overline{\mathcal{C}}=\oplus_{l \in \C} \mathcal{C}^l$ equipped with the bifunctor $\otimes$ and the trivial associativity constraint is a tensor category. \\
\paragraph{}For $V \in \mathcal{O}^l$, $W \in \mathcal{O}^{l'}$, $v \in V[\mu],\,w \in W[\nu]$ consider the formal composition 
$${\Phi}_{\lambda-\nu,k-l'}^v{\Phi}_{\lambda,k}^w:M_{\lambda,k} \to M_{\lambda-\mu-\nu,k-l-l'}{\otimes}V {\otimes}W$$
and the quantum correlation function is
$$\Psi_{\lambda,k}^{v,w}=<x_{\lambda-\mu-\nu,k-l-l'}^*,{\Phi}_{\lambda-\nu,k-l}^v{\Phi}_{\lambda,k}^w.x_{\lambda,k}>.$$
We now define the qKZ twist
 \begin{align*}
J_{V,W}(\lambda,k):\; V \otimes W & \to V \otimes W\\
v \otimes w &\mapsto \Psi_{\lambda,k}^{v,w}.
\end{align*}
\hbox to1em{\hfill}The following result is straightforward 
(it is a generalization of a similar result for finite dimensional Lie 
algebras, see \cite{EV1}):
\begin{theo} For $l \in \C$, the assignement 
\begin{align*}
F^l:\; {\mathcal O}^l &\to \mathcal{C}^l\\
V & \mapsto  (V,L_V(u,\hat\lambda)=R_{\C^n,V}(u,\lambda,k))
\end{align*}
and the identity at the level of morphisms is a well-defined functor. The system of functors $F^l$, $l \in \C$ gives rise to a tensor functor $\overline{F}:\mathcal{O}\to\overline{\mathcal{C}}$ with tensor structure
$$J_{V,W}(\hat\lambda): F(V) \otimes F(W) \to F(V \otimes W).$$
\end{theo}

In particular, if $V$ is a highest weight module over the quantum affine 
algebra, then $F(V)$ is a module with the same character over the 
elliptic quantum group. It would be interesting to describe these modules 
explicitly. Steps toward this goal are made in \cite{JKOS}. 
	
\section{Concluding remarks}
\paragraph{Generalization to arbitrary semisimple Lie algebra:} 
The results of this note generalize from $sl_n$ to any simple Lie algebra, 
with especially elegant formulas for classical Lie algebras. 
\paragraph{Link with $U_{p,q}(\hat{\mathfrak{sl}}_n)$:} Applying the Vertex-IRF transform to the category $\mathcal{C}^l$, as in \cite{ES} 
and references therein, one gets a functor from $\mathcal{C}^l$ to a certain category of representations of 
Belavin's quantum elliptic R-matrix with an arbitrary central charge, 
by matrix difference operators (where the matrices are infinite).
The defining RLL relations for such representations involves the Belavin 
R-matrix evaluated at two different values of the elliptic modulus $\tau$, 
as is the case for the algebra $U_{p,q}(\hat{\mathfrak{sl}}_n)$ (see \cite{FIJKMY}). Therefore we expect that using the methods of \cite{ES} and this note, 
one could establish an equivalence between the categories 
of highest weight modules 
over Belavin's and Felder's elliptic algebras.

\small{
}
\end{document}